\documentclass{amsart}
\usepackage{graphicx}
\usepackage[all]{xy}
\usepackage[ansinew]{inputenc}
\usepackage{amsmath}
\usepackage{amscd}
\usepackage{amssymb}
\usepackage{latexsym}
\usepackage{mathrsfs}
\usepackage{color}
\newtheorem{thm}{Theorem}

\theoremstyle{definition}

\newtheorem{lem}[thm]{Lemma}
\newtheorem{prop}[thm]{Proposition}
\newtheorem{defn}[thm]{Definition}

\newtheorem*{thmA}{Theorem A}
\newtheorem*{thmB}{Theorem B}
\newtheorem*{thmC}{Theorem C}
\newtheorem*{thmD}{Theorem D}
\newtheorem*{thmE}{Theorem E}

\numberwithin{equation}{section}
\newcommand{\N}{\mathbb{N}}
\newcommand{\Z}{\mathbb{Z}}
\newcommand{\Q}{\mathbb{Q}}
\newcommand{\R}{\mathbb{R}}

\newcommand{\Cal}{\mathcal}

\def \<{\langle}
\def \>{\rangle}

\def \((  {(\!(}
\def \)) {)\!)}

\begin{document}

\title[Interpreting the projective hierarchy]{Interpreting the projective hierarchy in expansions of the real line}%
\author{Philipp Hieronymi and Michael Tychonievich}
\address{University of Illinois at Urbana-Champaign\\
Department of Mathematics\\
1409 W. Green Street\\
Urbana, IL 61801\\
USA}
\email{P@hieronymi.de}

\thanks{A version of this paper will appear in the \emph{Proceedings of the American Mathematical Society}.}

\address{Department of Mathematics\\ The Ohio State University\\ 231 West 18th Avenue\\
Columbus, Ohio 43210\\ USA
}
\email{tycho@math.ohio-state.edu}


\date{\today}
%\dedicatory{}%
%\commby{}%

\maketitle

%\doublespacing

\begin{abstract} We give a criterion when an expansion of the ordered set of real numbers defines the image of $(\R,+,\cdot,\N)$ under a semialgebraic injection. In particular, we show that for a non-quadratic irrational number $\alpha$, the expansion of the ordered $\Q(\alpha)$-vector space of real numbers by $\N$ defines multiplication on $\R$.
\end{abstract}

\section{Introduction}

The main technical result is the following generalization of \cite[Theorem 1.1]{discrete}.

\begin{thmA}\label{mainthm} Let $\mathcal R$ be an expansion of the ordered set $(\R,<)$. If $\mathcal R$ defines an open interval $U\subseteq \R$, a closed and discrete set $D\subseteq \R_{\geq 0}$, and functions $f:D \to \R$ and $g: \R^3\times D\to D$  such that
\begin{itemize}
\item[(i)] $f(D)$ is dense in $U$,
\item[(ii)] for every $a,b \in U$ and $d,e \in D$ with $a<b$ and $e\leq d$
$$
\{ c\in \R : g(c,a,b,d)=e\}\cap \big(a,b\big) \hbox{ has non-empty interior,}
$$
\end{itemize}
then $\mathcal{R}$ defines every subset of $D^n$ and every open subset of $U^n$ for every $n\in \N$.
\end{thmA}
Hence a structure $\Cal R$ that satisfies the assumptions of Theorem A defines the image of $(\R,+,\cdot,\N)$ under a semialgebraic injection. The theory of such a structure is undecidable.\\

The significance of Theorem A comes from the fact that it is applicable outside the setting of expansions of the ordered \emph{field} of real numbers. In particular, Theorem A allows us to settle some open questions about expansions of the additive group of real numbers. The most significant is the following answer to questions raised by Chris Miller. A real number is called non-quadratic if it is not the solution to a quadratic equation with rational coefficients.

\begin{thmB} Let $\alpha \in \R$ be a non-quadratic irrational number and let $f: \R \to \R$ be the function that maps $x$ to $\alpha x$. Then $(\R,<,+,\N,f)$ defines multiplication on $\R$.
\end{thmB}

Note that the sets $\alpha \N$ and $\alpha^2 \N$ are definable in $(\R,<,+,\N,f)$. Since $\alpha$ is assumed to be non-quadratic, Theorem B follows immediately from the following stronger result.

\begin{thmC}
Let $\alpha,\beta,\gamma \in \R$ be such that $\alpha,\beta,\gamma$ are linearly independent over $\Q$. Then $(\R,<,+,\alpha\N,\beta \N,\gamma \N)$ defines multiplication on $\R$.
\end{thmC}

Theorems B and C are in stark contrast to results about $(\R,<,+,\N)$. This structure admits a quantifier elimination result and its theory is decidable (due independently to Weispfenning \cite{weis} and Miller \cite{ivp}). However, the structures in Theorem B and C define multiplication on $\R$ and hence every projective set. This is a much stronger property in respect to what kind of sets are definable than the property of just defining multiplication on $\N$. Consider $(\R,<,+,\N)^{\#}$, the expansion of $(\R, <,+,\N)$, by all nonempty subsets of all cartesian products $\N^k$, where $k$ ranges over all natural numbers. This structure defines multiplication on $\N$, is undecidable, but by Friedman and Miller \cite{sparse} every definable set is a union of an open set and a discrete set.\\

Theorem A can also be applied to expansions by analytic functions. Consider the structure $(\R,<,+,\sin)$. Marker and Steinhorn showed in unpublished work that every definable set in this expansion is a union of an open set and a discrete set\footnote{For expansions of $\R$, this property is equivalent to being locally o-minimal.}; for a proof of this result, see Toffalori and Vozoris \cite[Theorem 2.7]{sin}. Adding a predicate for $\N$ to this structure destroys this tameness completely.

\begin{thmD} The expansion $(\R,<,+,\sin,\N)$ defines multiplication on $\R$.
\end{thmD}

Several remarks about the proof of Theorem A and its connection to earlier work are in order. A stronger form of Theorem A was proved for expansions of the ordered fields in \cite{discrete}, but we will show that this result follows from Theorem A in this setting. The proof of Theorem A is not a minor modification of the proof of Theorem 1.1 in \cite{discrete}. The latter proof depends crucially on the definability of multiplication and does not transfer to our setting. In order to establish the conclusion of Theorem A, it is enough to show that every subset of all cartesian products $D^k$, where $k$ ranges over $\N$, is definable. For every $n\in\N$, boxes with endpoints in $f(D)$ form a countable basis of the induced topology on $U^n$. Hence we only need to show that the set of endpoints is definable. Similar to the earlier work in \cite{discrete} our approach relies on the idea of regarding the function $f$ as a definable approximation scheme. The key idea here is to realizes that for every $c\in U$, the set of best approximations of $c$ (see Definition 1) is definable and varies heavily in $c$. Given a subset $X\subseteq D^k$, in order to show that $X$ is definable, we describe a family of definable sets $(Y_a)$, depending on a parameter $a$, such that $Y_a$ is the image of a definable subset of the set of best approximations of $a$ under the function $g$. We then find a $c\in \R$ such that $Y_c = X$. The parameter $c$ will be constructed as a limit of a subset of $f(D)$. As in the earlier work, this construction rests on the topological completeness of $\R$.\\

The outline of the paper is as follows. In Section 2 some definitions and a few basic results are established. The proof of Theorem A and the details of how Theorem B to D follow from Theorem A can be found in Section 3. Section 4 contains a couple of remarks about the optimality of Theorem A.

\subsection*{Acknowledgements} The authors would like to thank Chris Miller for help in preparing this paper and the anonymous referee for very helpful remarks.

\subsection*{Notation} We will write $\N$ for $\{1,2,3,...\}.$

\section{Definitions}

In the following, let $\mathcal R$ be an expansion of the ordered set $(\R,<)$ such that $\mathcal R$ defines an open interval $U\subseteq \R$, a closed and discrete set $D\subseteq \R_{\geq 0}$ and a function $f:D \to U$ with $f(D)$ dense in $U$. In this section, we will introduce several definitions and establish some of their properties.

\begin{defn} Let $c\in U$. We say $d\in D$ is a \textbf{best approximation of $c$ from the left} if $f(d) < c$ and
$$
f(D_{<d}) \cap \big(f(d),c\big) = \emptyset.
$$
%
%$$
%\forall e \in D (e<d)\rightarrow (f(e)< f(d) \vee f(e) > c).
%$$
We write $L(c)$ for the set of best approximations of $c$ from the left. Similarly, we say $d \in D$ is a \textbf{best approximation of $c$ from the right} if $f(d)>c$ and
$$
f(D_{<d}) \cap \big(c,f(d)\big) = \emptyset.
$$
and write $R(c)$ for the set of best approximations of $c$ from the right.
\end{defn}

By density of $f(D)$ in $U$, $\sup f(L(c))=c=\inf f(R(c))$ for every $c\in U$.

\begin{defn} Let $d\in D$ and let $(L,R)$ be a pair of finite subsets of $D_{\leq d}$. We say $(L,R)$ is a \textbf{finite approximation} up to $d$ if
\begin{itemize}
\item[(1)] $f(d_l) < f(d_r)$, for $d_l \in L, d_r \in R$,
\item[(2)] $f$ is strictly increasing on $L$ and strictly decreasing on $R$,
\item[(3)]  for every $e \in D_{\leq d}$
%$$(\exists e_2 \in L \ e_2 \leq e_1 \wedge f(e_2) \geq f(e_1)) \vee (\exists e_2 \in R \ e_2 \leq e_1 \wedge f(e_2) \leq f(e_1)).
%$$
$$
\big(f(L\cap D_{\leq e})\cap [f(e),\infty)\big) \cup \big(f(R\cap D_{\leq e})\cap (-\infty, f(e)]\big) \neq \emptyset
$$

\end{itemize}
If $e \in D$ with $e>d$ and $(L',R')$ is an approximation up to $e$ such that $L\subseteq L'$ and $R\subseteq R'$, we say $(L',R')$ is an \textbf{extension} of $(L,R)$ up to $e$.
\end{defn}

\begin{lem}\label{approxcon} Let $(L,R)$ be a pair of finite subsets of $D$ and let $d \in D$. Then the following are equivalent
\begin{itemize}
\item [(i)] $(L,R)$ is a finite approximation up to $d$,
\item[(ii)] if $\max f(L) \leq c \leq \min f(R)$, then $$L(c) \cap D_{\leq d} = L \hbox{ and } R(c) \cap D_{\leq d} = R.$$
\end{itemize}
\end{lem}
\begin{proof} The proof of (ii)$\Rightarrow$(i) is immediate from the definition of a finite approximation. For the other implication, suppose that $(L,R)$ is a finite approximation.
Since $L$ and $R$ are finite, $\max f(L) < \min f(R)$ by (1). Let $c \in \R$ be such that $\max f(L) < c < \min f(R)$. Now (2) and (3) guarantee that $L(c) \cap D_{\leq d} = L$ and $R(c) \cap D_{\leq d} = R$. \end{proof}

\begin{lem} Let $(d_n)_{n\in \N}$ be a strictly increasing sequence of elements in $D$ and for every $n\in \N$, let $(L_{n},R_n)$ be a finite approximation up to $d_n$ such that both $L_n$ and $R_n$ are non-empty, and $(L_n,R_n)$ is an extension of $(L_{m},R_m)$ for every $m,n\in \N$ with $m<n$ . Then there is a unique $c\in U$ such that
\begin{itemize}
\item[(i)] $(L(c)\cap D_{\leq d_n},R(c)\cap D_{\leq d_n}) = (L_n,R_n)$ for all $n\in \N$, and
\item[(ii)] $\{c\} = \bigcap_{n\in \N} \big[\max f(L_n),\min f(R_n)\big]$.
\end{itemize}
In such a situation, we say $c$ is \textbf{approximated} by $\big( (L_n,R_n) \big)_{n\in\N}$.
\end{lem}
\begin{proof} Let $(d_n)_{n\in \N}$ and $\big( (L_n,R_n) \big)_{n\in\N}$ satisfy the assumptions. Since $L_m \subseteq L_n$ and $R_m \subseteq R_n$ for every $m,n \in \N$ with $m<n$, the sequence $(\max f(L_n))_{n\in \N}$ is increasing, the sequence $(\min f(R_n))_{n \in \N}$ is decreasing and $\max f(L_n) < \min f(R_n)$ for every $n\in \N$. Hence the set $\bigcap_{n\in \N} \big[\max f(L_n),\min f(R_n)\big]$ is a non-empty subset of $U$. Let $c$ be an element in this intersection. By Lemma \ref{approxcon}, $L_n = L(c) \cap D_{\leq d_n}$ and $R_n = R(c) \cap D_{\leq d_n}$ for every $n$. Hence
$$
\bigcup_{n\in \N} L_n = L(c) \hbox{ and } \bigcup_{n\in \N} R_n = R(c).
$$
Since $\sup f(L(c)) = \inf f(R(c))$, we get $\bigcap_{n\in \N} \big[\max f(L_n),\min f(R_n)\big]=\{c\}$.
\end{proof}

\begin{lem} Let $d_1,d_2 \in D$ with $d_1 \leq d_2$ and let $(L,R)$ be a finite approximation up to $d_1$. Then there is a unique extension $(L',R')$ of $(L,R)$ up to $d_2$ such that $L=L'$.
\end{lem}
\begin{proof} Let $c \in \R$ be such that $\max f(L)< c < \min f(R)$ and there is no $e \in D$ with $e \leq d_2$ and $\max f(L) \leq f(e) < c$. It is easy to see that $L(c) \cap D_{\leq d_1}=L$, $R(c) \cap D_{\leq d_1}=R$ and $L(c) \cap (d_1,d_2] = \emptyset$. Hence $(L,R(c) \cap D_{\leq d_2})$ extends $(L,R)$. The uniqueness follows from Lemma \ref{approxcon}.
\end{proof}

We call $(L',R')$ the \textbf{right extension} of $(L,R)$ up to $d_2$.

\begin{defn}Let $d_1,d_2 \in D$ with $d_1 < d_2$. We say $f$ \textbf{splits} between $d_1$ and $d_2$, if for every $e_1,e_2 \in D$ with $e_1,e_2\leq d_1$ and $f(e_1)<f(e_2)$, there is $e_3 \in D$ with $d_1 <e_3\leq d_2$ such that $f(e_1)<f(e_3)<f(e_2)$.
\end{defn}

By density of $f(D)$, for every $d_1\in D$ there is $d_2 \in D$ such that $d_1 < d_2$ and $f$ splits between $d_1$ and $d_2$. If $f$ splits between $d_1$ and $d_2$ and $d>d_2$, then $f$ splits between $d_1$ and $d$.

\section{Proofs}
\subsection*{Proof of Theorem A}
Let $\mathcal R$ be an expansion of the ordered set $(\R,<)$. Assume that $\Cal R$ defines an open interval $U\subseteq \R$, a closed and discrete set $D\subseteq \R_{\geq 0}$ and functions  $f:D \to R$ and $g: \R^3\times D\to D$  such that
\begin{itemize}
\item[(i)] $f(D)$ is dense in $U$,
\item[(ii)] for every $a,b \in \R$ and $d,e \in D$ with $a<b$ and $e\leq d$
$$
\{ c\in \R : g(c,a,b,d)=e\}\cap \big(a,b\big) \hbox{ has non-empty interior,}
$$
\end{itemize}
Let $W\subseteq U^m$ be open. We have to show that $W$ is definable. Since the lexicographic order topology on $\R^m$ has a countable basis and $f(D)$ is dense in $U$, we can find a
sequence  $((s_n,r_n))_{n\in \N}$ of elements in $D^{2m}$ such that
$$
\bigcup_{n \in \N} \big(f(s_{n,1}),f(r_{n,1})\big)\times ... \times \big(f(s_{n,m}),f(r_{n,m})\big) = W.
$$
It is only left to show that the range of these sequences is definable. We will show that for every $m\in \N$ every subset of $D^m$ is definable.\\

\noindent Let $A\subseteq D^m$ and let $\big((a_{i,1},...,a_{i,m})\big)_{i\in \N\cup \{0\}}$ be an enumeration of $A$. We will now show that $A$ is definable.\\

\noindent Let $c=(c_1,c_2,c_3) \in \R^3$. Let $\delta_{c}$ be the successor function on the well-ordered set $(L(c_1),<)$ and let $\delta_c^{i}$ denote the $i$-th iterate of $\delta_c$. For notational convenience, let $\delta_c^{0}$ be the identity function on $L(c_1)$. For $i=1,...,m$, define functions $l_{c,i} : L(c_1) \to \R$
by
$$
l_{c,i}(e) =\max f\big (L(c_3)\cap D_{\leq \delta_c^{i-1}(e)}\big)
$$
and $r_{c,i} : L(c_1) \to \R$ by
$$
r_{c,i}(e) =\min f\big (R(c_3)\cap D_{\leq \delta_c^{i-1}(e)}\big).
$$
Let $X_c$ be the set
$$
\{ e \in L(c_1) \ : \ L(c_2) \cap [e,\delta_c(e)] \neq \emptyset\}.
$$
Finally consider
$$
Y_c := \{(d_1,...,d_{m}) \in D^{m} : e \in X_c, d_i=g\big(c_3,l_{c,i}(e), r_{c,i}(e),\delta_c^{i-1}(e)\big) \}.
$$
Note that $Y_c$ is definable. The idea behind the above definition of $Y_c$ is to construct a tuple $c=(c_1,c_2,c_3) \in \R^3$ such that every $e\in L(c_1)$ picks $d_1 \in L(c_3)$ and $d_2 \in R(c_3)$ such that $g(c_3,d_1,d_2,e)$ has the desired value in $D$. The task of $L(c_2)$ is to code where a new tuple begins.

\begin{prop} There is $c\in \R^3$ such that $Y_c = A$.
\end{prop}
\begin{proof} We will construct a sequence $(d_n)_{n\in \N\cup \{0\}}$ of elements in $D$ and three strictly increasing sequences of finite approximations $\big((L_{i,n},R_{i,n})\big)_{n\in \N\cup\{0\}}$ for $i=1,2,3$, where $(L_{i,n},R_{i,n})$ is a finite approximation up to $d_n$ and
$(L_{i,n},R_{i,n})$ is an extension of $(L_{i,k},R_{i,k})$ for every $k,n\in \N\cup \{0\}$ with $k<n$. We will construct these four objects such that they have the following additional properties: for all $n\in \N\cup\{0\}$
\begin{itemize}
\item [(I$)_n$]$L_{1,n}=\{ d_k : k \in \N_{\leq n}\}$,
\item[(II$)_n$] $d_n \geq a_{i,j}$, for $i\leq n, j\leq m$,
\end{itemize}
and if $n\geq 1$ and  $n=ms+t$ for some $s,t \in \N\cup\{0\}$ with $1\leq t \leq m$,
\begin{itemize}
\item [(III$)_n$] $L_{2,n} \cap [d_{n-1},d_{n}] \neq \emptyset$ iff $t=1$,
 \item[(IV$)_n$] for all $x \in \big(\max f(L_{3,n}),\min f(R_{3,n})\big)$,
$$g(x,\max f(L_{3,n-1}), \min f(R_{3,n-1}),d_{n-1}) = a_{s,t}.$$
\end{itemize}

For $n=0$, let $d\in D$ be such that $d \geq a_{0,t}$ for every $t\leq m$. Let $d_0$ be the smallest element in $D$ larger than $d$ such that $f(d_0) < f(e)$ for all $e \in D$ with $e\leq d$. Take $a \in U$ such that $a>f(d_0)$, but $a < f(e)$ for all $e\in D$ with $e<d_0$.
 Set $L_{1,0} = L(a)\cap D_{\leq d_0}$ and $R_{1,0} = R(a) \cap D_{\leq d_0}$. By construction of $d_0$, $L_{1,0}=\{d_0\}$. Let $(L_{2,0},R_{2,0})$ be the right extension of $(\emptyset,\emptyset)$ up to $d_0$ and let $(L_{3,0},R_{3,0})$ be an arbitrary extension of $(\emptyset,\emptyset)$ up to $d_0$ such that both $L_{3,0}$ and $R_{3,0}$ are non-empty.\\

Suppose that $n\geq1$ and that $d_k$ and $\big((L_{i,k},R_{i,k})\big)$ are already constructed for $k < n$ and $i=1,2,3$ satisfying (I$)_k$, (II$)_k$ and, if $k\geq 1$, (III$)_k$ and (IV$)_k$. Let $s,t \in \N\cup \{0\}$ be such that $n=ms+t$ with $1\leq t \leq m$.
Take $d \in D$ such that
\begin{itemize}
\item[(A)] $d > d_{n-1}$ and $f$ splits between $d_{n-1}$ and $d$ and
\item[(B)] $d>\max \{ a_{i,j} : i \leq n, j=1,...,m\}$.
\end{itemize}
\noindent By (II$)_{n-1}$, $d_{n-1}>a_{s,t}$. By (ii), there are $$b_1,b_2 \in \big(\max f(L_{3,n-1}),\min f(R_{3,n-1})\big)$$ such that
 $g(x,\max f(L_{3,n-1}),\min f(R_{3,n-1}),d_{n-1}) = a_{s,t}$ for all $x\in \big(b_1,b_2)$. Take $b_3 \in \big(b_1,b_2\big)$ and $d'\in D$ larger than $d$ such that
 $$
 b_1 <  \max f(L(b_3)\cap D_{\leq d'}) < \min f(R(b_3)\cap D_{\leq d'}) < b_2.
 $$
 We will choose $(L_{3,n},R_{3,n})$ as an extension of $(L(b_3)\cap D_{\leq d'},R(b_3) \cap D_{\leq d'})$. Any such extension will satisfy (IV$)_n$. But first we have to choose $d_n$ and $(L_{i,n},R_{i,n})$ for $i=1,2$.\\

\noindent Let $(L_{1,n-1},R)$ be the right extension of $(L_{1,n-1},R_{1,n-1})$ up to $d'$. Let $d_{n}$ be the smallest element in $D$ larger than $d'$ such that
$$
\max f(L_{1,n-1}) < f(d_{n}) < \min f(R).
$$
Then $(L_{1,n-1}\cup \{d_{n}\},R)$ extends $(L_{1,n-1},R)$ and is a finite approximation up to $d_{n}$. Set $(L_{1,n},R_{1,n}):=(L_{1,n-1}\cup\{d_n\},R)$. Hence this extension satisfies (I$)_n$ and $d_n$ satisfies (II$)_n$.\\

\noindent
If $t\neq 1$, let $(L_{2,n},R_{2,n})$ be the right extension of $(L_{2,n-1},R_{2,n-1})$ up to $d_{n}$. \\
\noindent Now suppose that $t=1$. Since $f$ splits between $d_{n-1}$ and $d_{n}$, there is a minimal $e \in D$ with $d_{n-1}<e<d_{n}$ such that
$$\max f(L_{2,n-1}) <f(e) < \min f(R_{2,n-1}).$$
Then $(L_{2,n-1} \cup \{e\},R_{2,n-1})$ extends $(L_{2,n-1},R_{2,n-1})$ and is a finite approximation up to $e$. Let $(L_{2,n},R_{2,n})$ be an arbitrary extension of $(L_{2,n-1} \cup \{e\},R_{2,n-1})$ up to $d_{n}$. It follows directly that $(L_{2,n},R_{2,n})$ satisfies (III$)_n$.\\

\noindent Finally set $(L_{3,n},R_{3,n}):=(L(b_3)\cap D_{\leq d_{n}},R(b_3)\cap D_{\leq d_{n}})$. Since $d_n \geq d'$, $(L_{3,n},R_{3,n})$ is an extension of $(L(b_3)\cap D_{\leq d'},R(b_3) \cap D_{\leq d'})$. As mentioned above, such an extension satisfies (IV$)_n$.\\

For $i=1,2,3$, let $c_i \in U$ be the unique element that is approximated by $\big((L_{i,n},R_{i,n})\big)_{n\in \N\cup \{0\}}$ and let $c=(c_1,c_2,c_3)$. It is left to show that $Y_c=A$. Since $\big((L_{i,n},R_{i,n})\big)_{n\in \N\cup \{0\}}$ approximates $c_i$, we have by Lemma \ref{approxcon} that
$$L_{i,n} = L(c_{i}) \cap D_{\leq d_n} \hbox{ and } R_{i,n} = R(c_{i}) \cap D_{\leq d_n}.$$ By (I$)_n$,
$$
L(c_1) = \{ d_n : n \in \N\cup \{0\}\} \hbox{ and } \delta_c^i(d_n)=d_{n+i}.
$$
Let $n\in \N$ and $s,t \in \N$ such that $n=ms+t$ with $1\leq t \leq m$.
By (III$)_n$,
$$
L(c_2) \cap [d_{n-1},d_{n}] \neq \emptyset \hbox{ iff } t=1.
$$
Hence
$$
X_c = \{ d_{ms} : s \in \N\cup\{0\}\}.
$$
Since for every $k\in \N$
$$
\max f(L(c_3) \cap D_{\leq d_k}) = \max f(L_{3,k}) \hbox{ and } \min f(R(c_3) \cap D_{\leq d_k}) = \min f(R_{3,k}),
$$
we have for $i=1,...,m$
$$
l_{c,i}(d_{ms}) =\max f(L_{3,ms+i-1}) \hbox{ and } r_{c,i}(d_{ms}) =\min f(R_{3,ms+i-1}).
$$
By (IV$)_{ms+i}$,
\begin{align*}
g\big(c_3,l_{c,i}(d_{ms}),& r_{c,i}(d_{ms}),\delta_c^{i-1}(d_{ms})\big) = a_{s,i}.
\end{align*}
Thus $Y_c = A$.
\end{proof}

\subsection*{A first application of Theorem A} Let $\Cal R$ be an expansion of the ordered field of real numbers $(\R,<,+,\cdot)$. We will show that in this setting assumption (ii) in Theorem A follows from assumption (i).

\begin{thmE}  If $\Cal R$ defines a function $f:D \to U$ such that $D$ is closed and discrete, $U\subseteq \R$ is an open interval and $f(D)$ is dense in $U$, then
$\Cal R$ defines every open subset of $\R^n$ for every $n\in \N$.
\end{thmE}
\begin{proof} For every subinterval of $\R$, $\Cal R$ defines a bijection between this interval and $\R$. Hence we can assume that $U=(0,1)$. After replacing $D$ by
$$
\{ d \in D \ : \ \forall e \in D \ e<d\rightarrow f(e) \neq f(d)\},
$$
we can assume that $f$ is injective. We now construct a function $g: \R^3 \times D \to D$, definable in $\Cal R$, that satisfies condition (ii) of Theorem A.  First let $h: \R^2 \times D \to \R$ be given by
$$
h(a,b,e) := \left\{
                  \begin{array}{ll}
                    a + \frac{f(e)}{b - a} & \hbox{if } a<b, \\
                    a & \hbox{otherwise.}
                  \end{array}
                \right.
$$
Note that for fixed $a,b\in \R$ with $a<b$, the function
$
e \mapsto h(a,b,e)
$
is injective and its image is a subset of the interval $(a,b)$.
Then define $g : \R^3 \times D \to D$  such that $g(c,a,b,d)$ is the $e \in D_{\leq d}$  such that $|c-h(a,b,e)|$ is minimal. Since $D_{\leq d}$ is finite, $g$ satisfies condition (ii) in Theorem A.
\end{proof}

As a corollary, every expansion of the real field that satisfies the assumption of Theorem E defines $\N$. This was already shown in \cite[Theorem 1.1]{discrete}. Easy modifications of the above proof show that one needs only assume the definability of a homeomorphism between a bounded interval and an unbounded interval rather than the definability of multiplication.

\subsection*{Proof of Theorem C} We deduce Theorem C from Theorem A. For ease of notation, we assume that $\gamma=1$. So let $\alpha, \beta \in \R$ such that $1,\alpha, \beta$ are linearly independent over $\Q$. We will now show that $(\R,<,+,\N,\alpha\N,\beta \N)$ defines multiplication on $\R$.\\

Since for every $n\in \N$ the subgroups $n\alpha\N$ and $n\beta\N$ are definable, we can assume that $1<\alpha<\beta$. For every $a \in \R$, let $\lfloor a \rfloor_{\beta}$ be the largest element in $\beta\N$ smaller than $a$ if such an elements exists, and $0$ otherwise. Consider $f: \N \to (0,\beta)$ and $f_{\alpha} : \alpha \N \to (0,\beta)$, where
$$
f(n) := n - \lfloor n \rfloor_{\beta} \hbox{ and } f_{\alpha}(n\alpha) := n\alpha - \lfloor n \alpha \rfloor_{\beta}.
$$
Both $f$ and $f_{\alpha}$ are definable and the images of $f$ and $f_{\alpha}$ are dense in $(0,\beta)$ by Kronecker's Approximation Theorem (see \cite[Theorem 7.8]{apostol}). Define a function $h_1: \R_{>0} \times \N \times \N \to \alpha \N$ such that $h_1(u,d,e)$ is the minimum of the set of all $x\in \alpha\N$ such that
\begin{align*}
f_{\alpha}(x)\in \big(f(e),f(e)+u\big) \hbox{ and } f\big(\N \cap (0,d]\big) \cap \big(f(e),f_{\alpha}(x)\big] = \emptyset,
\end{align*}
in other words, given $(u,d,e)$, $h_1$ returns the smallest $x\in \alpha\N$ such that $f_{\alpha}(x)$ is larger than $f(e)$, but the difference is smaller than $u$ and there is no $y\in \N$ such that $y \leq d$ and $f(y)$ is between $f(e)$ and $f_{\alpha}(x)$. We will see that given $u$ and $d$, each $e\in \N_{\leq d}$ is uniquely determined by the distance between $f_{\alpha}(x)$ and $f(e)$.
Define a function $h_2 : \R_{>0} \times \N \times \N \to \Z + \alpha \Z+ \beta \Z$ by
$$
h_2(u,d,e) := f_{\alpha}(h_1(u,d,e))-f(e).
$$
For fixed $d \in \N$ and $u \in \R$, it follows directly from the linear independence of $1, \alpha, \beta$ over $\Q$ that the function that maps $e$ to $h_2(u,d,e)$ is injective on $\N \cap [0,d]$. Then let $g : \R^3 \times \N \to \N$ be defined such that if $b>a$, $g(c,a,b,d)$ is the $e \in \N \cap [0,d]$ such that
$$|(c-a)- h_2((b-a),d,e)|
$$
is minimal, and $g(c,a,b,d)=1$, if $b\leq a$. The function $g$ satisfies condition (ii) of Theorem A.\\

Hence by Theorem A, every subset of $\N^n$ is definable for every $n\in \N$. Let $\lfloor\cdot \rfloor : \R \to \Z$ be the usual floor function. Since $1<\alpha$, the floor function is injective on $\alpha \N$. Thus every subset of $\N^n \times (\alpha \N)^m$ is definable for every $n,m \in \N$. In particular, both multiplication on $\Z$ and scalar multiplication by $\alpha$ on $\Z$ is definable .\\

We will now show that multiplication by $\alpha^2$ is definable on $\Z$ as well. Since every subset of $\Z^3 \times (\alpha\Z)^2$ is definable, the map $\Z \times \Z \times \alpha \Z \to \Z \times \alpha \Z$ given by
$$
(n_1,n_2,\alpha m) \to (n_1n_2,\alpha n_1m)
$$
is definable. Thus the function $\Z \times (\Z +\alpha \Z) \to \Z +\alpha \Z$ that maps $(n_1,n_2+\alpha m)$ to $n_1n_2+\alpha n_1m$ is definable. Since this function is continuous and $\Z + \alpha \Z$ is a dense subset of $\R$, the function $\Z \times \R \to \R$ that maps $(n,a)$ to $na$ is definable. Fixing $\alpha^2$ in the second coordinate, we get that multiplication by $\alpha^2$ is definable on $\Z$.\\

In order to define multiplication on $\R$, it is enough to define multiplication on a dense subset of $\R$. We will show that multiplication is definable on the dense subset $\Z+\alpha \Z$. Since $(a+b)^2=a^2+2ab+b^2$, we just need to prove definability of the squaring function. Let $a \in \Z+\alpha \Z$.  Since there are unique $n,m \in \Z$ such that $a = n+\alpha m$ and multiplication by $\alpha$ is definable on $\Z$,  we can define the map
$a\mapsto (n,m)$. Since $a^2=n^2 + 2\alpha n m + \alpha^2 m^2$ and multiplication by $\alpha$ and $\alpha^{2}$ is definable of $\Z$, the squaring function is definable on $\Z + \alpha \Z$.

\subsection*{Proof of Theorem D} We have to show that $(\R,<,+,\sin,\N)$ defines multiplication on $\R$. As before, define a function $h_1: \R_{>0} \times \N \times \N \to \N$ such that $h_1(u,d,e)$ is the minimum of the set of all $x\in \N$ such that
$$
\sin(x) \in \big( \sin(e), \sin(e) + u\big) \hbox{ and } \sin(\N_{\leq d}) \cap \big(\sin(e),\sin(x)\big]=\emptyset.
$$
and define $h_2 : \R_{>0} \times \N \times \N \to \R$ by
$$
h_2(u,d,e) := \sin(h_1(u,d,e))-\sin(e).
$$
For fixed $d \in \N$ and $u \in \R$, it follows directly from the Lindemann-Weierstrass-Theorem (for a statement, see \cite[Theorem 1.4]{baker}) that the function that maps $e$ to $h_2(u,d,e)$ is injective on $\N_{\leq d}$. Then define $g : \R^3 \times \N \to \N$ to be such that for $b>a$,
$g(c,a,b,d)$ is the $e \in \N_{\leq d}$ such that $|(c-a)- h_2((b-a),d,e)|$ is minimal, and $g(c,a,b,d)=1$ for $b\leq a$. Now $g$ satisfies condition (ii) of Theorem A. The definability of multiplication on $\R$ follows as in the proof of Theorem C, using $2\pi\N$ instead of $\alpha\N$.

\section{A couple of remarks about optimality}
\subsection*{1} To our knowledge, there is no documented example of a structure that satisfies condition (i) of Theorem A, but neither condition (ii) nor the conclusion of the theorem. Hence it is not known whether condition (ii) follows from condition (i). While Theorem E shows that this implication holds in the setting of expansions of the ordered field of real numbers, we believe that it fails even for expansions of the additive group.

\subsection*{2} At the moment, we do not fully understand the complexity of the definable sets in $(\R,<,+,\alpha\N,\beta\N)$. In particular, we do not know whether the assumption or the conclusion of Theorem A holds for this structure, but we suppose that neither does. We have partial results that this structure defines complicated sets such as Cantor sets and infinitely branching trees. However since these are only partial results and the argument is of different nature than the ones presented in this paper, we decided not to include them here.

\end{document}